\documentclass[10pt]{amsart}

\usepackage[all]{xy}
\usepackage{amsmath,amssymb,amsfonts,amsthm,amsopn,dsfont}

\newtheorem{theorem}{Theorem}[section]

\newtheorem{example}[theorem]{Example}

\begin{document}

\title{A Two Dimensional Adler-Manin Trace and Bi-singular Operators}

\author {Farzad Fathizadeh, Masoud Khalkhali, Fabio Nicola, Luigi Rodino}

\date{}

\keywords{Adler-Manin trace, bi-singular operator, noncommutative residue}

\subjclass[2000]{58J42}

\address{F. Fathizadeh: 
Department of Mathematics and Statistics, 
York University, 
Toronto, Ontario, Canada, M3J 1P3,   
Current address: Department of Mathematics, 
Western University, 
London, Ontario, Canada, N6A 5B7 }
\email{ffathiz@uwo.ca}

\address{M. Khalkhali: Department of Mathematics, 
Western University, 
London, Ontario, Canada, N6A 5B7}
\email{masoud@uwo.ca}

\address{F. Nicola: Dipartimento di Scienze Matematiche, 
Politecnico di Torino, 
Corso Duca degli Abruzzi 24, 10129 Torino, Italy}
\email{fabio.nicola@polito.it}

\address{L. Rodino: Dipartimento di Matematica,  
 Universit\`{a} degli Studi di Torino, 
Via Carlo Alberto 10, 10123 Torino, Italy}
\email{luigi.rodino@unito.it}

\maketitle

\begin{abstract}

Motivated by the theory  of bi-singular pseudodifferential operators, 
we introduce a two dimensional version of the Adler-Manin trace. Our 
construction is rather general in the sense that it involves a twist afforded 
by an algebra automorphism. That is, starting from an algebra equipped 
with an automorphism, two twisted derivations, and a twisted invariant 
trace, we construct an algebra of formal twisted pseudodifferential symbols and 
define a noncommutative residue. Also, we provide related examples.

\end{abstract}

\tableofcontents

\section{Introduction}

The Adler-Manin trace was discovered for one dimensional 
pseudodifferential symbols \cite{adl, man}. It was generalized 
by Wodzicki to the algebra of pseudodifferential operators on a 
manifold of arbitrary dimension \cite{wod}. A vast generalization 
of Wodzicki's noncommutative residue was obtained by Connes 
and Moscovici for spectral triples in the context of the local index 
formula in noncommutative geometry \cite{conmos}. A pseudodifferential 
calculus was developed in \cite{con1} for $C^*$-dynamical systems. 
Noncommutative residues on the corresponding classical pseudodifferential 
operators for the canonical dynamical system defining noncommutative tori 
were studied in \cite{fatwon} (see also \cite{fatkha3}).   
Also, noncommutatvive residues for pseudodifferential operators on 
noncommutative tori with toroidal symbols are built in \cite{levjimpay}. 
Motivated by the notion of a twisted 
spectral triple \cite{conmos2}, a twisted version of the Adler-Manin trace was 
studied in \cite{fatkha}.

In this paper, motivated by the theory of bi-singular pseudodifferential 
operators \cite{duedue}, we find an algebraic setting for a noncommutative 
residue defined analytically in \cite{nic-rod}.  \par
Fredholm properties, index 
and residue were studied for an algebra of pseudodifferential operators, 
denoted by $HL(X_1\times X_2)$, on the product of two manifolds $X_1$ 
and $X_2$ in \cite{nic-rod}. 
In the the case $X_1=X_2=\mathbb{S}^1$, an operator in  $HL(X_1\times X_2)=HL(\mathbb{T})$ (with $\mathbb{T}=\mathbb{S}^1\times\mathbb{S}^1$) can be defined by a global symbol $a(x_1,x_2,\xi_1,\xi_2)$ as 
\begin{align}\label{uno-2013}
a(x_1,x_2,D_1,D_2)u&=\int e^{2\pi i (x_1\xi_1+x_2\xi_2)} a(x_1,x_2,\xi_1,\xi_2)\widehat{u}(\xi_1,\xi_2)\,d\xi_1\,d\xi_2\nonumber\\
&=\sum_{k_1,k_2\in\mathbb{Z}} e^{2\pi i (x_1k_1+x_ 2 k_2)} a(x_1,x_2,k_1,k_2)c_{k_1,k_2},
\end{align}
where $c_{k_1,k_2}=\iint_{\mathbb{T}} u(x_1,x_2) e^{-2\pi i (x_1k_1+x_2k_2)}\, dx_1 dx_2$ are the Fourier coefficients of $u\in C^\infty(\mathbb{T})$ (in the first line $u$ is regarded as a $1$-periodic function and its Fourier transform is a Dirac comb, which gives the discrete sum). The symbol $a$ is supposed to have an asymptotic expansion $a\sim\sum_{m=-\infty}^M\sum_{n=-\infty}^N \sigma_{m,n}(x_1,x_2,\xi_1,\xi_2)$, where $\sigma_{m,n}(x_1,x_2,\xi_1,\xi_2)$ is smooth in $\mathbb{T}\times(\mathbb{R}\setminus\{0\})^2$ and positively homogeneous of degree $m,n$ with respect to $\xi_1$ and $\xi_2$ respectively. Hence 
\[
\sigma_{m,n}(x_1,x_2,\xi_1,\xi_2)=\begin{cases} a_{m,n}^{(1,1)}(x_1,x_2)\xi_1^m\xi_2^n&{\rm if\ }\xi_1>0,\ \xi_2>0\\
a_{m,n}^{(1,2)}(x_1,x_2)\xi_1^m\xi_2^n&{\rm if\ }\xi_1>0,\ \xi_2<0\\
a_{m,n}^{(2,1)}(x_1,x_2)\xi_1^m\xi_2^n&{\rm if\ }\xi_1<0,\ \xi_2>0\\
a_{m,n}^{(2,2)}(x_1,x_2)\xi_1^m\xi_2^n&{\rm if\ }\xi_1<0,\ \xi_2<0\\
\end{cases}
\]
for suitable smooth functions $a^{(s,t)}_{m,n}$, ${s,t=1,2}$, on $\mathbb{T}$ (see  \cite[Section 2]{nic-rod} for more details). \par
Now, writing $I$ for the two sided ideal of the operators which are smoothing 
with respect to at least one of the variables on 
$\mathbb{T}$, we have an isomorphism
\begin{equation}\label{uno}
HL(\mathbb{T})/I \to \Psi (A,\delta_1,\delta_2),
\end{equation}
where $\Psi (A,\delta_1,\delta_2)$ is an algebra of formal pseudodifferential operators with coefficients in the algebra $A$, given by the direct sum of { $4$ copies} of the algebra $C^\infty(\mathbb{T})$, 
and $\delta_1,\delta_2$ are the standard partial derivatives. To be precise, each element in 
the right-hand side of \eqref{uno} is by definition a formal series 
\begin{equation}\label{due}
\sum_{m=-\infty}^M\sum_{n=-\infty}^N a_{m,n}\xi_1^m \xi_2^n,
\end{equation}
with $a_{m,n}=(a^{(s,t)}_{m,n})_{s,t=1,2}\in A$. We have 4 linearly 
independent traces  $\tau_{s, t}:A \to \mathbb{C}$ on the algebra $A$, defined by  
\[
\tau_{s,t} (a)=\iint_{\mathbb{T}} a^{(s,t)}(x_1,x_2)\,dx_1\, dx_2, \quad s,t=1,2,\quad a \in A. 
\]
Correspondingly, we will have 4 noncommutative residues on $\Psi (A,\delta_1,\delta_2)$, defined via $\tau_{s, t}, s, t=1,2,$ by considering the term $\tau_{s,t}(a_{-1,-1})$ for any $a \in A$.  We will prove these results in Section 4 for the {\it abstract} algebra $\Psi(A,\delta_1,\delta_2)$, generalizing therefore those obtained for bi-singular operators on the torus in \cite{nic-rod}.  \par
Bi-singular operators $P$ were originally studied, during the years 60's 
and 70's, under the form
\begin{eqnarray}\label{tre}
Pf(z_1,z_2) &=& b_0(z_1,z_2)f(z_1,z_2) \nonumber \\
&&+ \frac{1}{\pi i}\int_{\mathbb{S}^1}\frac{b_1(z_1,z_2,\zeta_1)}{\zeta_1-z_1}
f(\zeta_1,z_2)\,d\zeta_1 \nonumber \\
&& + \frac{1}{\pi i}\int_{\mathbb{S}^1}\frac{b_2(z_1,z_2,\zeta_2)}{\zeta_2-z_2}
f(z_1,\zeta_2)\,d\zeta_2 \nonumber \\ &&
-\frac{1}{\pi^2}\iint_{\mathbb{S}^1\times\mathbb{S}^1}
\frac{b_{1,2}(z_1,z_2,\zeta_1,\zeta_2)}{(z_1-\zeta_1)(z_2-\zeta_2)}
f(\zeta_1,\zeta_2)d\zeta_1\,d\zeta_2.
\end{eqnarray}
Here $z_1,z_2$ belong to $\mathbb{S}^1$, regarded as the unit circle in the complex plane, 
and we understood counter-clockwise integration in the principal value Cauchy sense. 
We assume the functions $b_0, b_1,$  $b_2, b_{1,2}$ are $C^\infty$ 
with respect to all of the variables. These operators, bounded on 
$L^2(\mathbb{T})$, can be regarded as operators of order zero in 
$HL(\mathbb{T})$, namely $a_{m,n}=0$ for $(m,n)\not=(0,0)$ in 
\eqref{due}. To respect the algebraic isomorphism in \eqref{uno}, one 
sets in \eqref{due} $a_{0,0}=(a^{(s,t)}_{0,0})_{s,t=1,2}$, where 
\begin{eqnarray}\label{quattro-2013}
a_{0,0}^{(s,t)}&=& b_0(z_1,z_2)-(-1)^s b_1(z_1,z_2,z_1) \nonumber \\ 
&& -(-1)^t b_2(z_1,z_2,z_2) +(-1)^{s+t} b_{1,2}(z_1,z_2,z_1,z_2), 
\end{eqnarray}
with $z_j= e^{2\pi i x_j},\ j=1,2$ (see Example \ref{esempio} below for computations). \par The literature in this connection is wide, involving problems of operator 
theory, harmonic analysis and several complex variables, see for example 
\cite{unodue,unosette,unouno}, coming from the school of A. Zygmund, and 
\cite{unotre,unoquattro,unosei,unocinque}, from the school of F.\ D.\ Gahov.

Starting from the 70's, bi-singular operators were reconsidered from the 
point of view of the theory of pseudodifferential operators, see 
\cite{dueuno} and \cite{duedue}, emphasizing the connection with the 
proof of the index theorem \cite{duetre}. Roughly, the algebra 
$HL(\mathbb{T})$ is obtained by composing $P$ in \eqref{tre} with 
partial derivatives on $\mathbb{T}$ and their inverses. This originates 
terms of positive and negative orders in \eqref{due}, and gives rise to 
non-trivial residues. Note that, arguing from a purely algebraic
point of view, Fredholm property and index have no 
meaning for $\Psi(A,\delta_1,\delta_2)$, nevertheless a trace $\tau$ on
$\Psi(A,\delta_1,\delta_2)$ in \eqref{uno} can be obviously extended to 
$HL(\mathbb{T})$ by setting $\tau=0$ on the ideal  $I$.

A natural question is whether a related class of pseudodifferential 
operators exists, providing in \eqref{uno} an isomorphism with 
$\Psi(A,\delta_1,\delta_2)$ where simply $A=C^\infty(\mathbb{T})$. 
To this end, we may address to the Toeplitz operators in a quadrant 
of \cite{treuno,tredue,tretre,trequattro}.

This paper is organized as follows. In Section \ref{AdlManTrace},  
we recall some standard  constructions and the 
Adler-Manin trace in the one dimensional case. In Section \ref{twistedncres}, we 
recall from \cite{fatkha} the extension of the Adler-Manin trace to a twisted set-up, where 
the twist is afforded by an algebra automorphism. This work was motivated by the notion 
of twisted spectral triples, introduced recently by Connes and Moscovici \cite{conmos2},  
which appear naturally in the study of type III examples of foliation algebras. They have 
shown that given a twisted spectral triple, the Dixmier trace induces a twisted trace 
on the base algebra, where lack of any non-trivial trace is the characteristic property 
of such algebras. Also,  we briefly recall from \cite{mos} the construction of spectral triples 
twisted by scaling automorphisms, and the existence of a tracial noncommutative residue on 
the corresponding algebra of twisted pseudodifferential operators.  In Section \ref{tdnr}, we 
construct a two dimensional version of the algebra of twisted formal pseudodifferential 
symbols and introduce a noncommutative residue on this algebra. That is, starting from an 
algebra equipped with an automorphism and two twisted derivations, we construct an algebra of 
formal twisted pseudodifferential symbols. Then, we introduce a noncommutative residue on 
this algebra via an invariant twisted trace on the base algebra, and show that it is a trace 
functional. In Section \ref{Examples}, as an example of our two dimensional analogue of the 
Adler-Manin noncommutative residue, we shall consider the noncommutative residue of  
bi-singular operators.

\section{The Adler-Manin Trace} \label{AdlManTrace}

Let $A$ be a unital associative algebra over $\mathbb{C}$ which 
is not necessarily commutative. Recall that a \emph{derivation} on 
$A$ is a linear map $\delta : A \to A$ such that 
\[
\delta (ab)= a \delta (b) + \delta (a) b, \qquad a,b \in A, 
\] 
Given a pair $(A, \delta)$ as above, the algebra 
of \emph{formal differential symbols} $D(A, \delta)$ \cite{kas1,man} is, 
by definition, the algebra generated by $A$ and a symbol $\xi$ subject to 
the relations
\begin{equation} \label{eq:relation}
\xi a = a \xi + \delta (a), \qquad a \in A. 
\end{equation}
For every element $D$ of $D(A, \delta)$ there is a unique expression of the form
\[
\sum_{i=0}^{N} a_i \xi^i, 
\]
with $N \geq 0, a_i \in A$.  We think of $D$ as a differential 
operator of order at most $N$. Using \eqref{eq:relation}, one can prove 
by induction that
\begin{equation} \label{eq:relation2}
\xi^n a = \sum_{j=0}^{n} \binom{n}{j} \delta^j(a) \xi^{n-j},  
\qquad a \in A, \qquad n \geq 0.
\end{equation}
Using \eqref{eq:relation2}, we obtain the following multiplication 
formula in $D(A, \delta)$:
\[
\Big ( \sum_{i=0}^{M}a_i \xi^i \Big ) \Big ( \sum_{j=0}^{N} b_j \xi^j \Big ) 
= \sum_{n=0}^{M+N} \Big ( \sum_{i,j,k} \binom{i}{k} a_i \delta^{k}(b_j) \Big ) \xi^n, 
\]
where the internal summation is over all $0 \leq k  \leq i \leq M$, and 
$0 \leq j \leq N$ such that $i+j-k=n$.

By formally inverting $\xi$ in $D(A, \delta)$  and completing the resulting 
algebra, one obtains the algebra of \emph{formal pseudodifferential symbols} 
of $(A, \delta)$ \cite{adl,man,kas1}, denoted by $\Psi (A,\delta)$. More 
precisely, it is defined as follows. Elements of $\Psi (A,\delta)$ consist 
of formal sums
\[
D=\sum_{i=- \infty}^{N} a_i \xi^i, 
\]
with $a_i \in A$, and $N \in \mathbb{Z}$. Its multiplication is defined 
by extending \eqref{eq:relation2} to all $n \in \mathbb{Z}$ as follows. For any 
$a \in A$ and $n \in \mathbb{Z}$ one defines 
\begin{equation}
\xi^n a = \sum_{j=0}^{\infty} \binom{n}{j} \delta^j(a) \xi^{n-j}. \nonumber
\end{equation}
Here the binomial coefficient $\binom{n}{j}$ for $n \in \mathbb{Z}$ and 
non-negative $j \in \mathbb{Z}$, is defined by $\binom{n}{j}:=
\frac{n(n-1) \cdots (n-j+1)}{j!}$. Notice that for $n < 0$, we have an 
infinite formal sum. It follows that for general  $D_1 = \sum_{i=- \infty}^{M}
a_i \xi^i$ and $D_2 = \sum_{j=- \infty}^{N} b_j \xi^j$ in $\Psi(A, \delta)$, the 
multiplication is given by
\[
D_1 D_2 = \sum_{n=-\infty}^{M+N} \Big ( \sum_{i,j,k} 
\binom{i}{k} a_i \delta^{k}(b_j) \Big ) \xi^n, 
\]
where the internal summation is over all integers $i \leq M$, $ j \leq N$, 
and $ k \geq 0 $ such that $i+j-k =n$.

Now consider a $\delta$-\emph{invariant trace}  $\tau : A \to \mathbb{C}$. By definition, 
$\tau$ is a linear functional such that 
\[
\tau (ab) = \tau (ba), \qquad \tau (\delta (a))=0, \qquad  a,b \in A.
\]  
The \emph{Adler-Manin noncommutative residue}
\cite{adl,man,kas1} is the linear functional  $\text{res} : \Psi (A, \delta) \to \mathbb{C}$ 
defined by
\[
\text{res}  \Big ( \sum_{i=-\infty}^{N}a_i \xi^i \Big ) = \tau (a_{-1}).
\]
One checks that $\text{res}$ is a trace, i.e.
\[
\text{res} ([D_1,D_2] )= \text{res}(D_1D_2 - D_2 D_1)= 0, 
\qquad D_1, D_2 \in \Psi (A, \delta).
\]
Equivalently, one shows that the map
$\text{res} : \Psi (A, \delta) \to A/([A,A]+ \text{im} \, \delta)$
\[
D \mapsto a_{-1} \qquad \text{mod} \qquad [A,A] + \text{im} \, \delta
\]
is a trace on $\Psi(A,\delta)$ with values in $A/([A,A] + \text{im} \, \delta)$.

A relevant example is when $A=C^{\infty}(\mathbb{S}^1)$, the algebra of 
smooth functions on the circle  $\mathbb{S}^1= \mathbb{R} / \mathbb{Z}$,
with 
\[
\tau (f) = \int_{0}^{1} f(x)\,dx, \qquad  \delta(f)=f', \qquad f \in C^{\infty}(\mathbb{S}^1).
\] 
In this case the noncommutative residue coincides with the Wodzicki residue on 
the algebra of classical pseudodifferential operators on the circle, see Section 
\ref{Examples} below.

\section{Twisted Symbols and  Noncommutative Residues} 
\label{twistedncres}

Motivated by the notion of a twisted spectral triple \cite{conmos},  the 
Adler-Manin trace was extended to a twisted set-up in \cite{fatkha}. That is, 
an algebra of twisted formal pseudodifferential symbols was 
defined, where the twist is afforded by an automorphism of the 
base algebra, and a noncommutative residue was introduced. 
In this section we briefly recall this construction and highlight the 
tracial property of the corresponding noncommutative residue. 
We also recall some facts about noncommutative geometric spaces 
\cite{acbook}, namely spectral triples and the corresponding 
algebras of pseudodifferential operators, which admit noncommutative 
residues under some mild conditions \cite{conmos}. We also briefly recall 
the twisted spectral spectral triples obtained from scaling automorphisms 
of a spectral triple in \cite{mos}, and the corresponding algebras of twisted 
pseudodifferential operators and noncommutative residues.

\subsection{Adler-Manin trace in a twisted set-up}

We assume that $A$ is a unital associative algebra over 
$\mathbb{C}$, and $\sigma: A \to A$ is an algebra automorphism. 
A linear map $\delta: A \to A$ is a $\sigma$-{\it derivation} or a {\it twisted 
derivation} if 
\[
\delta(ab)=\delta(a) b+ \sigma(a) \delta(b), \qquad a,b \in A.
\]
Given such a triple $(A, \sigma, \delta)$, the algebra of {\it formal twisted 
differential symbols} $D(A, \sigma, \delta)$ is defined as follows. Its elements 
are polynomials in $\xi$ with coefficients from $A$, which are of the form
\[
\sum_{i=0}^N a_i \xi^i, \qquad N \geq 0, \qquad a_i \in A.
\]
The multiplication of   $D(A, \sigma, \delta)$ is defined by the relations 
\[
\xi a = \sigma(a) \xi + \delta(a), \qquad a \in A.
\]
By induction, it follows that
\begin{equation} \label{ntwisted}
\xi^n a = \sum_{i=0}^{n} P_{i, \, n}(\sigma , \delta)(a) \xi^{i}, \qquad n >0, \qquad a \in A, 
\end{equation}
where $P_{i,\, n}(\sigma , \delta) : A \to A$ is the noncommutative 
polynomial in $\sigma$ and $\delta$ with $\binom{n}{i}$ terms of 
total degree $n$ such that the degree of $\sigma$ is $i$. For example 
\[
P_{3,4}(\sigma, \delta) = \delta \sigma^3 + \sigma \delta \sigma^2 + 
\sigma^2 \delta \sigma + \sigma^3 \delta. 
\]
By formally inverting $\xi$ and using induction, for $n=-1$ and $a \in A,$ we have 
\[ 
\xi^{-1} a= \sum_{i=0}^{N} (-1)^i \sigma^{-1} (\delta 
\sigma^{-1})^{i}(a) \xi^{-1-i} + (-1)^{N+1} \xi^{-1} (\delta 
\sigma^{-1})^{N+1}(a)\xi^{-1-N}.
\]
This suggests setting 
\begin{equation} \label{eq:twistedrelation-1}
\xi^{-1} a= \sum_{i=0}^{\infty} (-1)^i \sigma^{-1} (\delta 
\sigma^{-1})^{i}(a) \xi^{-1-i}. \nonumber
\end{equation}
Therefore for  $n > 0, a \in A,$ we set
\begin{eqnarray} \label{-ntwisted}
&& \xi^{-n}a =  \nonumber \\ 
&& \qquad \sum_{i_1=0}^{\infty} \cdots \sum_{i_{n}=0}^{\infty} 
(-1)^{i_1  + \cdots + i_{n}} \sigma^{-1} (\delta \sigma^{-1})^{i_{n}} 
\cdots \sigma^{-1} (\delta \sigma^{-1})^{i_1} (a) \xi^{-n - i_1 - 
\cdots - i_{n}}.  
\end{eqnarray}

Thus, the elements of the algebra of {\it formal twisted 
pseudodifferential symbols} \cite{fatkha}, denoted by $\Psi(A, \sigma, \delta)$, 
are of the form 
\[
\sum_{-\infty}^N a_i \xi^i, \qquad N \in \mathbb{Z}, \qquad a_i \in A,   
\]
and the multiplication in this algebra is defined as follows. 
For any 
\[D_1= \sum_{n=- \infty}^{N} a_n \xi^n, \quad  
D_2 = \sum_{m=- \infty}^{M} b_m \xi^m \quad \in \quad \Psi(A, \sigma, \delta),  
\] 
we have 
\begin{eqnarray} \label{twistedmultiplication}
D_1 D_2 &=&
 \sum_{m= - \infty}^{M} \sum_{n<0} \sum_{i \geq 0} (-1)^{| i |} a_n 
\sigma^{-1} (\delta \sigma^{-1})^{i_{-n}} \cdots \sigma^{-1} 
(\delta \sigma^{-1})^{i_1} (b_m) \xi^{m+n- |i|} {} \nonumber \\
&& + \sum_{m= - \infty}^{M} \sum_{n=0}^{N} \sum_{j=0}^{n} 
a_n P_{j,\, n} (\sigma, \delta) (b_m) \xi^{m+j}, \nonumber
\end{eqnarray} 
where for $n<0$, $i=(i_1, \dots , i_{-n})$ is an $n$-tuple of integers 
and $|i|= i_1 + \cdots + i_{-n}$.

One of the main results proved in \cite{fatkha} is that given a twisted 
$\delta$-invariant trace $\tau : A \to \mathbb{C}$, the linear functional 
$\text{res}_\sigma : \Psi(A, \sigma, \delta) \to \mathbb{C}$ defined by 
\[
\text{res}_\sigma \Big ( \sum_{n=- \infty}^{N} a_n \xi^n \Big ) = \tau(a_{-1}) 
\]
is a trace. Here, by {\it twisted trace} we mean that 
\[
\tau (ab) = \tau (\sigma (b)a), \qquad  a,b \in A.
\]
Also, $\tau$ is said to be {\it $\delta$-invariant} if $\tau \circ \delta =0$.

\subsection{Spectral triples twisted by scaling automorphisms} 
\label{sttsa}

Geometric spaces are described by {\it spectral triples} $(\mathcal{A}, \mathcal{H}, D)$ 
in noncommutative geometry \cite{acbook, conmos}. Here, $\mathcal{A}$ is a 
$*$-algebra represented by bounded operators on a Hilbert space $\mathcal{H}$, 
and $D$ is an unbounded selfadjoint operator acting in $\mathcal{H}$, which 
interacts with $\mathcal{A}$ in a bounded fashion. That is, for any $a \in \mathcal{A}$, 
the commutator $[D, a]=Da - aD$ extends by continuity to a bounded linear operator on 
$\mathcal{H}$. A local index formula is proved for spectral triples in \cite{conmos}. If a 
spectral triple is $n$-{\it summable}, the Dixmier trace $\textnormal{Tr}_\omega$ induces 
a trace on the base algebra $\mathcal{A}$. This trace functional is given by 
\[
a \mapsto \textnormal{Tr}_\omega(a |D|^{-n}), \qquad a \in \mathcal{A}. 
\]
Also, it is shown in \cite{conmos} that if the spectral triple has a {\it simple 
dimension spectrum}, the analogue of Wodzicki's residue defines a 
trace on the corresponding algebra  $\Psi(\mathcal{A}, \mathcal{H}, D)$ of 
pseudodifferential operators. This trace is defined by 
\[
P \mapsto \textnormal{Res}_{z=0} \textnormal{Trace} (P |D|^{-2z}), 
\qquad P \in \Psi(\mathcal{A}, \mathcal{H}, D).
\]

Since the existence of a non-trivial trace is the characteristic of type II 
situations in Murray-von Neumann classification of rings of operators, 
in order to incorporate type III examples, the notion of twisted spectral 
triples was introduced recently by Connes and Moscovici \cite{conmos2}.  
In fact, this notion arises naturally in the study of type III examples of 
foliation algebras, and also in noncommutative conformal geometry 
\cite{conmos2, contre, fatkha1, conmos3, fatkha2}. For twisted spectral triples, 
the ordinary commutators $[D, a]$ are not necessarily bounded operators, 
however, there exists an algebra automorphism $\sigma: \mathcal{A} \to \mathcal{A}$ 
such that the twisted commutators $[D, a]_\sigma=Da - \sigma(a) D$ extend to bounded 
operators for all $a \in \mathcal{A}$.

Twisted spectral triples also arise naturally in conformal geometry of Riemannian 
manifolds \cite{conmos2, mos}.  Let $(M, g)$ be a connected compact 
Riemannian spin manifold of dimension $n$ and $D = D_g$ 
be the associated Dirac operator acting on the Hilbert space of  $L^2$-spinors 
$\mathcal{H}= L^2(M, S^g)$. Let $SCO(M, [g])$ denote the Lie group of 
diffeomorphisms of $M$ that preserve the conformal structure $[g]$
(consisting of all Riemannian metrics that are conformally equivalent to $g$), 
the orientation, and the spin structure. Also, let $G = SCO(M, [g])_0$ denote the 
connected component of the identity. In \cite{mos}, using a suitable automorphism 
of the crossed product algebra $C^{\infty}(M) \rtimes G$, a twisted spectral triple of
the form $(C^{\infty}(M) \rtimes G, L^2(M, S^g), D)$ is constructed.

Similarly, by endowing $\mathbb{R}^n$ with the Euclidean metric, and 
considering the group $G$ of conformal transformations of $\mathbb{R}^n$, 
a twisted spectral triple  is constructed over the crossed product algebra 
$C_c^{\infty}(\mathbb{R}^n) \rtimes G$ \cite{mos}. An abstract formulation 
of this class of twisted spectral triples leads to the idea of twisting an ordinary 
spectral triple by its \emph{scaling automorphisms}.  A local index formula 
is proved in \cite{mos} for these examples under a particular invariance 
property under the twist, which enforces the Selberg principle for classical examples.

The set of scaling automorphisms of a spectral triple 
$(\mathcal{A}, \mathcal{H}, D)$, denoted by $\textnormal{Sim}(\mathcal{A}, 
\mathcal{H}, D)$,
consists of all unitary operators $U$ on $\mathcal{H}$ such that
\[
U \mathcal{A} U^*= \mathcal{A} \qquad \textnormal{and} \qquad UDU^*=\mu(U)D \qquad
\textnormal{for some} \qquad\mu(U) >0. \nonumber 
\]
One can see that $\textnormal{Sim}(\mathcal{A}, \mathcal{H}, D)$ is a group 
and the map $\mu :\textnormal{Sim}(\mathcal{A}, \mathcal{H}, D) \to (0, \infty)$ 
is a character. Since the group $\textnormal{Sim}(\mathcal{A}, \mathcal{H}, D)$ 
acts by conjugation on $\mathcal{A}$, by fixing a subgroup
$G \subset \textnormal{Sim}(\mathcal{A}, \mathcal{H}, D)$, one can form the crossed  
product algebra  $\mathcal{A}_G = \mathcal{A} \rtimes G$. By definition, as a vector 
space, $\mathcal{A}_G$ is equal to $\mathcal{A} \otimes \mathbb{C} G$ whose elements 
are finite sums of the form 
\[ 
\sum a_U U, \qquad U \in G, \qquad a_U \in \mathcal{A}.  
\]
The multiplication in this algebra is defined by 
\[
U a = (UaU^*) U, \qquad U \in G, \qquad a \in \mathcal{A}. 
\]
It is shown in \cite{mos} that the formula
\[
\sigma(aU)= \mu(U)^{-1}aU, \qquad  a\in \mathcal{A}, \qquad U \in G, \nonumber 
\]
defines an automorphism of $\mathcal{A}_G$, and $(\mathcal{A}_G, \mathcal{H}, D)$ 
is a twisted spectral triple. For the twisted commutators, one has
\[ [D,aU]_{\sigma}=DaU - \sigma(aU) D =[D, a ]U,  \nonumber \]
which are bounded operators for all $a \in \mathcal{A}, U \in G$.

For this class of twisted spectral triples, the group  
$\textnormal{Sim}(\mathcal{A}, \mathcal{H}, D)$ acts by conjugation 
on the algebra of pseudodifferential operators $\Psi(\mathcal{A},
\mathcal{H}, D)$ of the base spectral triple. In \cite{mos} the crossed product algebra 
\[
\Psi(\mathcal{A}\rtimes G, \mathcal{H}, D) := \Psi(\mathcal{A}, \mathcal{H}, D) \rtimes G
\] 
is considered, and it is shown that the residue functional given by
\[
P \mapsto \textnormal{Res}_{z=0} \textnormal{Trace} (P |D|^{-2z}), \qquad 
P \in \Psi(\mathcal{A}\rtimes G, \mathcal{H}, D), 
\]
defines a trace functional, under a similar assumption to the untwisted case, namely the 
\emph{extended simple dimension spectrum hypothesis}.

\section{Two Dimensional Noncommutative Residue}\label{tdnr}

In this section, motivated by the theory of bi-singular pseudodifferential 
operators \cite{duedue}, we find an algebraic setting for a noncommutative 
residue defined analytically in \cite{nic-rod}. The setting is a higher-dimensional 
version of the one for the Adler-Manin trace, explained in Section \ref{AdlManTrace}. 
Our construction is rather general in the sense that, similar to the 
construction described in Section \ref{twistedncres}, it involves a twist afforded 
by an algebra automorphism.

As above, we assume that  $A$ is a unital associative algebra and 
$\sigma: A \to A$ is an algebra automorphism. We consider two 
$\sigma$-derivations $\delta_1, \delta_2 : A \to A$ such that 
$\delta_1 \delta_2 = \delta_2 \delta_1, \delta_1 \sigma = \sigma \delta_1, 
\delta_2 \sigma = \sigma \delta_2$. First, we associate to this data an algebra 
$D(A, \sigma, \delta_1, \delta_2)$ of twisted differential symbols whose elements 
are polynomials in $\xi_1, \xi_2$ with coefficients in $A$, which are of the form
\[
\sum_{m=0}^M \sum_{n=0}^{N} a_{m, n} \xi_1^m \xi_2^n, \qquad M, N >0, 
\qquad a_{m,n} \in A.
\]
The multiplication in $D(A, \sigma, \delta_1, \delta_2)$ is essentially defined 
by the relations 
\[
\xi_i a = \sigma(a) \xi_i + \delta_i(a), \qquad i=1,2, \qquad a \in A.
\]
Since $\delta_i$ and $\sigma$ commute, by induction, for any $n >0$ we have 
\[
\xi_i^n a = \sum_{j=0}^n \binom{n}{j} \sigma^{n-j} \delta_i^{j} (a) \xi_i^{n-j}, 
\qquad i=1,2, \qquad a \in A.
\]

We also define an algebra of formal twisted pseudodifferential symbols 
$\Psi(A, \sigma, \delta_1, \delta_2)$ (in short $\Psi(A, \delta_1, \delta_2)$ if $\sigma=id$) formally inverting each $\xi_i, i=1,2$. 
The elements of this algebra are  of the form
\[ 
\sum_{m=-\infty}^M \sum_{n=-\infty}^{N}  a_{m,n} \xi_1^m \xi_2^n, \qquad 
M, N \in \mathbb{Z}, \qquad a_{m,n} \in A. 
\]
The multiplication of this algebra is defined by the following relations:
\[
\xi_1 \xi_2 = \xi_2 \xi_1,
\]
\[  
\xi_i \xi_i^{-1} = \xi_i^{-1} \xi_i = 1, \qquad i=1,2, 
\]
\[ 
\xi_i^n a = \sum_{j=0}^\infty \binom{n}{j} \sigma^{n-j} \delta_i^{j} (a) \xi_i^{n-j}, 
\qquad i=1,2, \qquad n \in \mathbb{Z}, \qquad a \in A.
 \]
We note that  these relations follow from \eqref{ntwisted} and \eqref{-ntwisted} 
under the assumption that the automorphism and the derivations commute with 
each other, which is assumed in our construction as mentioned above.   
Therefore, for  any 
\[
D_1=\sum_{m=-\infty}^M \sum_{n=-\infty}^N a_{m,n} \xi_1^m \xi_2^n, \quad 
D_2=\sum_{p=-\infty}^{M'} \sum_{q=-\infty}^{N'} b_{p,q} \xi_1^p \xi_2^q \quad \in \quad
\Psi(A, \sigma, \delta_1, \delta_2),
\] we 
have: 
\begin{eqnarray}
&& D_1D_2 = \nonumber \\ 
&& \qquad \sum \sum_{j_1, j_2 \geq 0} 
\binom{m}{j_1} \binom{n}{j_2} a_{m,n} 
\sigma^{m+n-j_1-j_2} \delta_1^{j_1} \delta_2^{j_2} (b_{p,q})
\xi_1^{m+p-j_1} \xi_2^{n+q-j_2}, \nonumber 
\end{eqnarray}
where the first summation is over all integers $m, n, p, q$ such that 
$m \leq M, n \leq N, p \leq M', q \leq N'.$

In the following, we define a noncommutative residue 
$\textnormal{Res}: \Psi(A, \sigma, \delta_1, \delta_2) \to \mathbb{C}$ 
and prove that it is a trace functional. For the definition, we start 
from a twisted trace $\tau: A \to \mathbb{C}$, which satisfies some 
invariance properties stated in the hypotheses of Theorem \ref{nctrace}.

\newtheorem{first}{Definition}[section]
\begin{first}
For any  
\[
D= \sum_{m=-\infty}^M\sum_{n=-\infty}^N a_{m,n} \xi_1^m \xi_2^n \quad  \in \quad
\Psi(A, \sigma, \delta_1, \delta_2), 
\]
we define its noncommutative residue by  
\[
\text{Res}(D) = \tau(a_{-1,-1}).
\]
\end{first}

In order to investigate the tracial property of $\textnormal{Res}$, which is 
carried out in Theorem \ref{nctrace}, first, we need to prove a lemma:

\newtheorem{intparts}[first]{Lemma}
\begin{intparts} \label{intparts}
Let $\delta:A \to A $ be a $\sigma$-derivation such that 
$\delta \sigma = \sigma \delta$. If $\tau: A \to \mathbb{C}$ is 
a linear functional such that $\tau \circ \delta= 0$, then for any $a, b \in A$ 
and non-negative integer $i,$ we have
\[ 
\tau( \delta^i(a)b) = (-1)^i \tau (\sigma^i(a) \delta^i(b)).
\]
\begin{proof}
For any $a, b \in A$, we have
\begin{eqnarray}
\tau(\delta(ab)) &=&  \tau(\delta(a)b+\sigma(a) \delta(b)) \nonumber \\
&=&\tau(\delta(a)b)+ \tau(\sigma(a) \delta(b)) \nonumber \\
&=&0. \nonumber
\end{eqnarray}
Therefore
\[ 
\tau (\delta(a)b) = - \tau( \sigma(a) \delta(b)), \qquad a, b \in A.
\]
So the statement holds for $i=1$.  By induction, for any non-negative integer $i$ we have: 
\begin{eqnarray}
\tau(\delta^i(a)b) &=& \tau \big (\delta^{i-1} (\delta(a)) b \big ) \nonumber \\
&=&  (-1)^{i-1} \tau \big (\sigma^{i-1}(\delta(a)) \delta^{i-1}(b) \big ) \nonumber \\ 
&=& (-1)^{i-1} \tau \big (\delta (\sigma^{i-1}(a)) \delta^{i-1}(b) \big ) \nonumber \\
&=& (-1)^i \tau (\sigma^i(a) \delta^i(b)). \nonumber
\end{eqnarray}
\end{proof}
\end{intparts}

Now, we show that under suitable invariance properties on the twisted 
trace $\tau: A \to \mathbb{C}$, the noncommutative residue defined above 
gives a trace on the corresponding algebra of two dimensional formal twisted 
pseudodifferential symbols.

\newtheorem{nctrace}[first]{Theorem}
\begin{nctrace} \label{nctrace}
Let $\tau: A \to \mathbb{C}$ be a $\sigma^2$-trace such that 
$\tau \circ \delta_i = 0$ for $i=1,2$, and $\tau \circ \sigma=\tau$. 
Then the noncommutative residue $\text{Res}$ gives a trace functional 
on  $ \Psi(A, \sigma, \delta_1, \delta_2)$.

\begin{proof}

Since $\textnormal{Res}$ is clearly a linear functional, in order to prove that it 
is a trace, it suffices to show that for 
any $a, b \in A$ and $m, n, p, q \in \mathbb{Z}$, we have: 
\begin{equation}
\textnormal{Res} (a \xi_1^m \xi_2^n b \xi_1^{p} \xi_2^{q}) = 
\textnormal{Res}(b \xi_1^{p} \xi_2^{q} a \xi_1^m \xi_2^n). \nonumber
\end{equation}
First, we observe that 
\begin{eqnarray}
a \xi_1^m \xi_2^n b \xi_1^{p} \xi_2^{q} &=& 
a \xi_1^m \sum_{j=0}^\infty \binom{n}{j} \sigma^{n-j} \delta_2^j(b)  
\xi_1^{p} \xi_2^{n-j+q} \nonumber \\
&=&  a  \sum_{j=0}^\infty \binom{n}{j} \xi_1^m \sigma^{n-j} \delta_2^j(b)  
\xi_1^{p} \xi_2^{n-j+q} \nonumber \\
&=& a \sum_{j=0}^\infty \binom{n}{j}  \sum_{i=0}^\infty \binom{m}{i}  
\sigma^{m-i}\delta_1^i \sigma^{n-j} \delta_2^j(b) \xi_1^{m-i+p}  \xi_2^{n-j+q} \nonumber \\
&=&  \sum_{i, j=0}^\infty  \binom{m}{i}  \binom{n}{j} a \,
\sigma^{m-i+n-j}\delta_1^i \delta_2^j(b) \xi_1^{m+p-i}  \xi_2^{n+q-j}. \nonumber
\end{eqnarray}
Therefore
\begin{equation} \label{lside}
\textnormal{Res} (a \xi_1^m \xi_2^n b \xi_1^{p} \xi_2^{q})  = 
\tau \Big( \sum \binom{m}{i}  \binom{n}{j} a \,
\sigma^{m-i+n-j}\delta_1^i \delta_2^j(b) \xi_1^{m+p-i}  \xi_2^{n+q-j} \Big),
\end{equation}
where the summation is over all non-negative integers $i, j$ such that  
$m+p-i= -1$ and $n+q-j =-1$.

Similarly, we have:
\[
b \xi_1^{p} \xi_2^{q} a \xi_1^m \xi_2^n= 
\sum_{i, j=0}^\infty  \binom{p}{i}\binom{q}{j} b \,\sigma^{p-i+q-j}
\delta_1^i \delta_2^j(a) \xi_1^{m+p-i} \xi_2^{n+q-j}. 
\]
Therefore
\begin{eqnarray} \label{rside}
\textnormal{Res}(b \xi_1^{p} \xi_2^{q} a \xi_1^m \xi_2^n) &=& 
\tau \Big (    \sum  \binom{p}{i}\binom{q}{j} b 
\,\sigma^{p-i+q-j}\delta_1^i \delta_2^j(a)  \Big ), \nonumber
\end{eqnarray}
where the summation is again over all non-negative integers 
$i, j$ such that $m+p-i= -1$ and $n+q-j =-1$. Since $\tau$ is a 
$\sigma^2$-trace, for each term in \eqref{rside} we have 
\begin{eqnarray}
\tau \Big ( \binom{p}{i}\binom{q}{j} b \,\sigma^{p-i+q-j}\delta_1^i \delta_2^j(a) \Big ) &=&  
\binom{p}{i}\binom{q}{j} \tau ( \sigma^{p-i+q-j+2}\delta_1^i \delta_2^j(a) b) \nonumber \\
&=& \binom{p}{i}\binom{q}{j} \tau ( \delta_1^i \delta_2^j \sigma^{p-i+q-j+2} (a) b), \nonumber
\end{eqnarray}
which, using Lemma \ref{intparts}, is equal to 
\[
 \binom{p}{i}\binom{q}{j} (-1)^{i+j}\tau ( \sigma^{p+q+2} (a) \delta_1^{i} \delta_2^{j}(b) ).
\]
Using $\tau \circ \sigma = \tau$, the latter is equal to 
\[
\binom{p}{i}\binom{q}{j} (-1)^{i+j}\tau (  a \sigma^{-p-q-2} \delta_1^{i} \delta_2^{j}(b) ).
\]
Since $m+p-i=n+q-j=-1$ in the above sum, we have $-p-q-2=m+n-i-j.$ 
Thus we have 
\begin{equation} \label{rside}
\textnormal{Res}(b \xi_1^{p} \xi_2^{q} a \xi_1^m \xi_2^n) = 
\tau \Big (    \sum  \binom{p}{i}\binom{q}{j} (-1)^{i+j}\tau 
(  a \sigma^{m+n-i-j} \delta_1^{i} \delta_2^{j}(b) )  \Big ),
\end{equation}
where the summation is again over all non-negative integers 
$i, j$ such that $m+p-i= -1$ and $n+q-j =-1$.

By comparing  \eqref{lside} and \eqref{rside}, in order to prove that they 
are identical, we show that 
\[ 
\binom{m}{i}\binom{n}{j}  =(-1)^{i+j}\binom{p}{i}\binom{q}{j}, 
\]
for any fixed non-negative integers $i, j$ such that $m+p-i= -1$ and 
$n+q-j =-1$. This  can be  proved by showing that 
$(-1)^i\binom{p}{i}=\binom{m}{i}$ and similarly 
$ (-1)^j\binom{q}{j}= \binom{n}{j}$ as follows. We have:

\begin{eqnarray}
(-1)^i \binom{p}{i} &=& (-1)^i \binom{-m+i-1}{i} \nonumber \\
&=& (-1)^i \frac{(-m+i-1)(-m+i-1-1) \cdots (-m+i-1 -i+1)}{i !} \nonumber \\
&=& \frac{(m-i+1)(m-i+2) \cdots (m)}{i !} \nonumber \\
&=& \binom{m}{i}. \nonumber
\end{eqnarray}
\end{proof}
\end{nctrace}

\begin{example}\label{exa31}\rm
Let $A= C^\infty (\mathbb{T})$, $\sigma=id$, $\delta_1=\frac{\partial}{\partial \theta_1}$, 
$\delta_2=\frac{\partial}{\partial \theta_2}$, and
\[
\tau(f)=\int_0^1\int_0^1 f(\theta_1,\theta_2)\,d\theta_1\,d\theta_2, \qquad f \in A.
\]
 Then
 \[
 {\rm Res}\Big( \sum_{m=-\infty}^{M}\sum_{n=-\infty}^{N} 
 a_{m,n}\xi_1^m\xi_2^n  \Big)=\int_0^1\int_0^1 a_{-1,-1}(\theta_1,\theta_2)\,d\theta_1\,d\theta_2.
 \]
\end{example}

\section{Singular and Bi-Singular Operators}\label{Examples}

We add here some detail concerning the bi-singular operators mentioned in 
the introduction. To be definite, we consider first the one dimensional case.

\subsection{Singular integral operators}

Classical pseudodifferential operators \cite{agranovich,ruzhansky} of order 
$N\in\mathbb{Z}$ on the circle $\mathbb{S}^1=\mathbb{R}/\mathbb{Z}$  can be defined by 
a global symbol $a(x,\xi)$, $x\in \mathbb{R}/\mathbb{Z}$, $\xi\in\mathbb{R}$, 
with an asymptotic expansion in homogeneous terms. Modulo regularizing 
operators, i.e. mappings $\mathcal{D}'(\mathbb{S}^1)\to C^\infty(\mathbb{S}^1)$, 
this algebra can therefore be identified with that of formal series
\begin{equation}\label{otto}
\sigma=\sum_{n=-\infty}^N a_{n,1}(x)\xi_+^n+a_{n,2}(x)\xi_-^n,
\end{equation}
where $\xi_+=\max\{\xi,0\}$, $\xi_-=\min\{\xi,0\}$.\par
This corresponds to the Adler-Manin algebra 
\[
\Psi(A,\delta)=\Big\{\sum_{n=-\infty}^N a_n \xi^n;  \quad N \in \mathbb{Z},  \quad a_n \in A  \Big\}, 
\] 
with $A= C^\infty (\mathbb{S}^1)\oplus C^\infty(\mathbb{S}^1)$, $\delta=\frac{d}{dx}$, 
$a_n=(a_{n,1},a_{n,2})$.  There are two linearly independent traces on $A$, namely
\[
\tau_{s}(a)=\int_0^1 a^{(s)}(x)\,dx,\qquad s=1,2,\qquad a=(a^{(1)},a^{(2)})\in A. 
\]
Therefore, there are two noncommutative residues
\[
{\rm Res}_{s}\Big( \sum_{n=-\infty}^N a_n \xi^n\Big)=\tau_{s}(a_{-1})=\int_0^1 
a_{n,s}(x)\,dx,\qquad s=1, 2,
 \]
which are known as Wodzicki's
residues.\par
Introducing on $\Psi(A,\delta)$ the Toeplitz projection
\[
P_+\sigma=\sum_{n=-\infty}^N a_{n,1}(x)\xi_+^n,
\]
i.e. setting $a_{n,2}=0$ for all $n$ in \eqref{otto}, we may define the quotient 
subalgebra of the asymptotic expansions of the classical Toeplitz operators, 
cf. the example at the end of Section 2, where $A$ is simply given by 
$C^\infty(\mathbb{S}^1)$.

\subsection{Bi-singular integral operators} 
Consider the algebra $HL(\mathbb{T})$ of bi-singular
integral operators on
$\mathbb{T}=\mathbb{S}^1\times\mathbb{S}^1$, which are
defined by a global symbol  
\[a(x_1,x_2,\xi_1,\xi_2),\]
as in \eqref{uno-2013}, 
where $x_1,x_2\in \mathbb{R}/\mathbb{Z}$,
$\xi_1,\xi_2\in\mathbb{R}$, with an asymptotic expansion in
terms which are positively homogeneous with respect to $\xi_1$ and
$\xi_2$ separately. As explained in the introduction, modulo sums of operators regularizing
in $x_1$ or $x_2$, they are identified with
formal series of the form 
\[
\sigma=\sum_{m=-\infty}^M\sum_{n=-\infty}^N \sigma_{m,n},
\]
where  each $\sigma_{m,n}$ is given by an expansion of the form 
\[ 
\sigma_{m,n}=\sum_{s,t=1,2} a_{m,n}^{(s,t)}(x_1,x_2)\xi_{1,s}^m\xi_{2,t}^n,
\]
with $\xi_{j,1}=\max\{\xi_j,0\}$, $\xi_{j,2}=\min\{\xi_j,0\}$, $j=1,2$.\par
Therefore, the above operators correspond to the algebra
\begin{eqnarray}
&& \Psi(A, \delta_1,\delta_2)=\Psi(A, \sigma, \delta_1,\delta_2) = \nonumber \\ 
&& \qquad \Big\{ 
\sum_{m=-\infty}^M\sum_{n=-\infty}^N a_{m,n}\xi_1^m\xi_2^n; 
\quad M, N \in \mathbb{Z}, \quad a_{m,n}=(a_{m,n}^{(s,t)})_{s,t=1,2}\in A   
\Big\}, \nonumber
\end{eqnarray}
considered in Section \ref{tdnr}, with 
\[ 
A=C^\infty( \mathbb{T})^4, \qquad \sigma=id, \qquad \delta_1
=\frac{\partial}{\partial x_1},  \qquad \delta_2=\frac{\partial}{\partial x_2}.
\] 
Now we have four linearly independent traces on $A$, namely
\[
\tau_{s,t}(a)=\int_0^1\int_0^1 a^{(s,t)}(x_1,x_2)\,dx_1\, dx_2,\qquad 
a=(a^{(s,t)})\in A,\qquad s,t=1,2. 
\]
Therefore, in this case, we have four noncommutative residues on 
$\Psi(A, \sigma, \delta_1,\delta_2) $, which are 
given by 
\[
{\rm Res}_{s,t}\Big( \sum_{m=-\infty}^M \sum_{n=-\infty}^N 
a_{m,n} \xi_1^j\xi_2^k \Big) = \tau_{s,t}(a_{-1,-1}),\qquad s,t=1,2.
\]
They correspond to those found in \cite{nic-rod} for general
bi-singular operators on the product of two manifolds.\par
Finally, we note that the projection 
\[
P_{++}\sigma=\sum_{m=-\infty}^M\sum_{n=-\infty}^N 
a_{m,n}^{(1,1)}(x_1,x_2)\xi_{1,1}^m\xi_{2,1}^n
\]
gives the algebra of the asymptotic expansions for Toeplitz operators 
in the quarter plane, corresponding to Example \ref{exa31}.
\begin{example}\label{esempio}\par\rm
We finally verify that the operator $P$ in \eqref{tre} falls in the class $HL(\mathbb{T})$ and prove the formula \eqref{quattro-2013} for its symbol. \par
To this end, in \eqref{tre} we perform a $0$th order Taylor expansion of $b_1,$  $b_2,\  b_{1,2}$ at $\zeta_1=z_1$, $\zeta_2=z_2$ and $\zeta_1=z_1$, $\zeta_2=z_2$ respectively. Ignoring operators whose integrals kernels are smooth with respect to a couple of variables, which belong to the ideal $I$, we get 
\begin{eqnarray*}\label{tre-bis}
Pf(z_1,z_2) &=& b_0(z_1,z_2)f(z_1,z_2) + \frac{b_1(z_1,z_2,z_1)}{\pi i}\int_{\mathbb{S}^1}\frac{f(\zeta_1,z_2)}{\zeta_1-z_1}
\,d\zeta_1\nonumber \\
&& +\frac{b_2(z_1,z_2,z_2)}{\pi i}\int_{\mathbb{S}^1}\frac{f(z_1,\zeta_2)}{\zeta_2-z_2}
\,d\zeta_2 \\
&&
-\frac{b_{1,2}(z_1,z_2,z_1,z_2)}{\pi^2}\iint_{\mathbb{S}^1\times\mathbb{S}^1}
\frac{f(\zeta_1,\zeta_2)}{(z_1-\zeta_1)(z_2-\zeta_2)}
d\zeta_1\,d\zeta_2\ {\rm mod}\ I.
\end{eqnarray*}
The formula \eqref{quattro-2013} then follows by expressing the Hilbert transform on a circe as a pseudodifferential operator:
\[
\frac{1}{\pi i}\int_{\mathbb{S}^1} \frac{u(\zeta)}{\zeta-z}\,d\zeta=\sum_{k\geq 0} e^{2\pi i xk} a(x,k)c_{k}-\sum_{k<0} e^{2\pi i xk} a(x,k)c_{k}=a(x,D) u+Ru,
\]
where $c_{k}=\int_{\mathbb{S}^1} u(x) e^{-2\pi i xk}\, dx$, $R$ is a regularizing operator, and $a(x,\xi)=(1-\chi(\xi)){\rm sign}\, \xi$, $\chi$ being a smooth cut off function, $\chi(\xi)=1$ in a neighborhood of the origin.

\end{example}


\begin{thebibliography}{9}


\bibitem{adl} 
M. Adler,  
{\it On a trace functional for formal pseudo differential operators and 
the symplectic structure of the Korteweg-de\thinspace Vries type equations,}   
Invent. Math.  50  (1978/79), no. 3, 219--248.


\bibitem{agranovich} 
M.S. Agranovich, 
{\it Spectral properties of elliptic pseudodifferential operators on a closed curve, }  
(Russian) 
Funktsional Anal. i Prilo\v{z}en 13 (1979), no. 4, 54--56. 
(English translation in Functional Analysis and Its applications. 13, p. 279--281).


\bibitem{duetre} 
M.\ F. Atiyah and I.\ M. Singer, 
{\it The index of elliptic operators,}  
Ann. of Math.  87 (1968), 484--530.


\bibitem{tretre} 
L.\ A. Coburn, R.\ G. Douglas and I.\ Singer, 
{\it An index theorem for Wiener-Hopf operators on the discrete quarter-plane. }
J. Differential Geometry 6 (1972), 587-593.


\bibitem{con1} A. Connes,  
{\it $C^*$-alg\`ebres et g\'eom\'etrie diff\'erentielle}, 
C.R. Acad. Sc. Paris, t.~290, S\'erie A, 599--604, 1980.


\bibitem{acbook} 
A. Connes, 
{\it Noncommutative Geometry.} Academic Press, 1994.


\bibitem{conmos} 
A. Connes and H.  Moscovici,  
{\it The local index formula in noncommutative geometry,} 
Geom. Funct. Anal. 5 (1995), no. 2, 174--243.


\bibitem{conmos2}
A. Connes and H.  Moscovici,
{\it Type III and spectral triples,} 
Traces in number theory, geometry and quantum fields, 57--71, 
Aspects Math., E38, Friedr. Vieweg, Wiesbaden, 2008. 


\bibitem{conmos3} A. Connes, H. Moscovici, 
{\it Modular curvature for noncommutative two-tori,} 
arXiv:1110.3500.


\bibitem{contre} A. Connes, P. Tretkoff, 
{\it The Gauss-Bonnet theorem for the noncommutative two torus},
Noncommutative geometry, arithmetic, and related topics, 141-158, 
Johns Hopkins Univ. Press, Baltimore, MD, 2011.


\bibitem{treuno} 
R.\ G. Douglas and R. Howe, 
{\it On the $C^\ast$-algebras of Toeplitz operators on the quarter plane,} 
Trans. Amer. Math. Soc. 158 (1971), 203--217.


\bibitem{tredue}  
R.\ G. Douglas, 
{\it On the invertibility of a class of Toeplitz operators on the quarter-plane,} 
Indiana Univ. Math. J. 21 (1971/1972), 1031--1035.


\bibitem{dueuno} 
R.\ V. Dudu\v{c}ava, 
{\it On bisingular integral operators with discontinuous coefficients,} 
Mat. Sbornik 101 (1976); transl. Math. USSR Sbornik 30 (1976), 515--537.


\bibitem{fatkha} 
F. Fathizadeh and M. Khalkhali,  
{\it The algebra of formal twisted pseudodifferential  symbols and a noncommutative residue,} 
Letters  in  Mathematical Physics (2010), 94:41--61.


\bibitem{fatkha1} F. Fathizadeh, M. Khalkhali, 
{\it The Gauss-Bonnet theorem for noncommutative two tori 
with a general conformal structure}, 
J. Noncommut. Geom. 6 (2012), no. 3, 457--480.


\bibitem{fatkha2} F. Fathizadeh, M. Khalkhali, 
{\it Scalar curvature for the noncommutative two torus,} 
To appear in the Journal of Noncommutative Geometry, 
arXiv:1110.3511.


\bibitem{fatkha3} F. Fathizadeh, M. Khalkhali, 
{Scalar Curvature for Noncommutative Four-Tori,} 
arXiv:1301.6135. 



\bibitem{fatwon} F. Fathizadeh, M. W. Wong, 
{\it Noncommutative residues for pseudo-differential operators on the 
noncommutative two-torus}, 
Journal of Pseudo-Differential Operators and Applications,  2(3) (2011), 289--302.


\bibitem{unouno} C. Fefferman, 
{\it Estimates for double Hilbert transforms,} 
Studia Math. 44 (1972), 1--15.


\bibitem{unocinque} 
V.A. Kaki\v{c}ev. 
{\it Regularization and pseudoregularization of complete bisingular operators 
with Cauchy kernels and normal characteristic part,} 
Dokl. Akad. Nauk SSSR 229 (1976), 19--22; transl. Soviet Math Dokl. 17 (1976), 952--956.


\bibitem{kas1} 
Ch. Kassel,  
{\it Le r\'{e}sidu non commutatif} (d'apr$\grave{\textrm{e}}$s M. Wodzicki), 
(French) [The noncommutative residue (after M. Wodzicki)] 
S\'{e}minaire Bourbaki, Vol. 1988/89.  Ast\'{e}risque  No. 177-178  (1989), Exp. No. 708, 199--229.


\bibitem{levjimpay}
C. L\'evy, C. Neira Jim\'enez, S. Paycha, 
{\it The canonical trace and the noncommutative residue on the noncommutative torus,} arXiv:1303.0241. 


\bibitem{man} 
Ju. I. Manin, 
{\it Algebraic aspects of nonlinear differential equations. 
(Russian)  Current problems in mathematics, }
Vol. 11,  pp. 5--152. (errata insert) Akad. Nauk SSSR Vsesojuz. Inst. 
Nau\v cn. i Tehn. Informacii, Moscow, 1978.


\bibitem{mos}
H. Moscovici, 
{\it Local index formula and twisted spectral triples,} 
Quanta of maths, 465--500, Clay Math. Proc., 11, Amer. Math. Soc., Providence, RI, 2010.

\bibitem{nic-rod} 
F. Nicola and L. Rodino, 
{\it Residue index for bisingular operators,}  
in ``$C^\ast$ algebras and Elliptic Theory", Trends in Mathematics (2006), 187--202.


\bibitem{unoquattro} 
V.\ S. Pilidi, 
{\it The index of bisingular operators,} 
Funkcional Anal. i Prilo\v{z}en 7 (1973), 93--94; transl. Functional Anal. Appl. 7 (1973), 337--338.


\bibitem{unosei} 
V.\ S. Pilidi, L.\ I. Sazonov, 
{\it A priori estimates for characteristic bisingular integral operators,} 
Dokl. Akad. Nauk SSSR 217 (1974), 285--287; transl. Soviet Math. Dokl. 15 (1974), 1064--1067.


\bibitem{duedue} 
L. Rodino, 
{\it A class of pseudodifferential operators on the product of two manifolds and applications,} 
Ann. Scuola Norm. Sup. Pisa 2 ser. IV (1975), 287--302.


\bibitem{ruzhansky} 
M. Ruzhansky and V. Turunen, 
{\it Pseudodifferential operators and symmetries,}  
Pseudodifferential Operators, Theory and Applications, Birkh\"auser, Basel, 2010.


\bibitem{unotre} 
I.\ B. Simonenko, 
{\it On the question of solvability of bisingular and polysingular equations,}  
Funkcional Anal. i Prilo\v{z}en 5 (1971), 93--94; transl. Functional Anal. Appl. 5 (1971), 81--83.


\bibitem{trequattro} 
G. Strang, 
{\it Toeplitz operators in a quarter-plane,} 
Bull. Amer. Math. Soc. 76 (1970), 1303--1307.


\bibitem{unosette} 
T. Walsh, 
{\it On the existence of double singular integrals for kernels without smoothness,} 
Proc. Amer. Math. Soc. 28 (1971), 439--445.


\bibitem{wod}
M. Wodzicki, 
{\it Noncommutative residue. I. Fundamentals, }
$K$-theory, arithmetic and geometry (Moscow, 1984--1986),  320--399,
Lecture Notes in Math., 1289, Springer, Berlin, 1987.


\bibitem{unodue} A. Zygmund, 
{\it On the boundary values of functions of several complex variables I,}  
Fund. Math. 36 (1949), 207--235.


\end{thebibliography}
\end{document}